\documentclass[11pt]{article}

\usepackage[margin=1in]{geometry}
\usepackage[T1]{fontenc}
\usepackage[utf8]{inputenc}
\usepackage{amsmath,amssymb,amsthm,mathtools}
\usepackage{bm}
\usepackage{graphicx}
\usepackage{booktabs}
\usepackage{authblk}
\usepackage[numbers,sort&compress]{natbib}
\usepackage[colorlinks=true,linkcolor=blue,citecolor=blue,urlcolor=blue]{hyperref}
\usepackage{cleveref}

\theoremstyle{plain}
\newtheorem{theorem}{Theorem}[section]

\newtheorem{corollary}[theorem]{Corollary}
\theoremstyle{definition}

\theoremstyle{remark}
\newtheorem{remark}[theorem]{Remark}

\input{macros}

\title{\bf Minimax rates for pointwise CATE estimation:\\
the role of design smoothness}

\author[1]{Yuxuan Hou}
\author[2]{Zhenghan Li}
\affil[1]{Qiuzhen College, Tsinghua University}
\affil[2]{School of Statistics, University of Minnesota}

\date{\today}

\begin{document}
\maketitle

\begin{abstract}
We study the role of design smoothness in the minimax estimation of conditional
average treatment effects at a fixed covariate point. For design densities bounded
above and away from zero, without further regularity, we derive upper and lower
risk bounds with an arbitrarily small gap in the convergence-rate exponent. Under
uniform design, we obtain matching upper and lower bounds at interior points in
the regime of low nuisance smoothness. These results show that rough covariate
designs can lead to slower rates than uniform designs, even under the same
smoothness assumptions on the treatment effect and nuisance functions,
highlighting the role of the covariate distribution in pointwise causal estimation.
\end{abstract}


\medskip
\noindent\textbf{MSC2020 subject classifications:} Primary 62G08, 62C20;
secondary 62G20.


\section{Introduction}
\label{sec:intro}
The conditional average treatment effect (CATE) describes how the average effect of a
treatment varies with observed covariates. We study its estimation at a fixed covariate
value, with accuracy measured by mean absolute error.
Existing approaches include matching
\citep{rosenbaum1983central,heckman1998matching,abadie2006large,abadie2011bias,gao2020minimax,lin2023estimation}
and doubly robust or orthogonal estimation
\citep{robins1994estimation,bang2005doubly,lee2017doubly,nie2021quasi,kennedy2023towards}.
Our focus is on how the covariate distribution affects the accuracy attainable by any estimator.

For pointwise CATE estimation over H\"older classes, \citet{kennedy2024minimax}
established minimax lower bounds and constructed a higher-order local polynomial
R-learner that attains these bounds under suitable conditions on the design density.
Their lower-bound construction uses a design density that vanishes on parts of the
covariate domain. This leaves open whether the same lower bound holds under a uniform
density. A second question is what rate is attainable when the density is only bounded
above and away from zero, without the regularity needed for accurate density estimation.

Under uniform design, we establish the lower bound of \citet{kennedy2024minimax}
at interior covariate points in the regime of low nuisance smoothness. The nuisance
functions are the conditional treatment probability and the control outcome regression.
We also derive a matching upper bound under conditions on the nuisance estimators and
the accuracy of design density estimation, and verify these conditions with high
probability under uniform design. For rough design, where the density is bounded above
and away from zero without smoothness restrictions, we give upper and lower risk bounds
with an arbitrarily small gap in the convergence-rate exponent. These bounds identify
low-smoothness regimes in which a rough design slows pointwise estimation.

Our conditional upper bound follows the higher-order local polynomial R-learner
framework of \citet{kennedy2024minimax}. The construction combines local polynomial
estimation with a second-order projection correction based on a design density estimated
from independent training data. For rough design, we use weighted contrasts of nearby
observations to cancel low-order polynomials without estimating the density.

Section~\ref{sec:setup} introduces the model and the two design settings.
We present the rough-design bounds in Section~\ref{sec:rough} and the
uniform-design results in Section~\ref{sec:smooth}.
Section~\ref{sec:discussion} discusses extensions.

\section{Setup}
\label{sec:setup}

\subsection{Model and parameter classes}
\label{sec:model}
We estimate the conditional mean difference $\tau(x)=\mu_1(x)-\mu_0(x)$, where
$\mu_a(x)=\E(Y\mid X=x,A=a)$ for $a\in\{0,1\}$.
Under consistency, conditional exchangeability, and the overlap condition
below, this conditional mean difference identifies the CATE.

\paragraph{H\"older classes.}
For $\zeta>0$ let $\underline\zeta=\lceil\zeta\rceil-1$ and $\theta_\zeta=\zeta-\underline\zeta\in(0,1]$. For
$g:[0,1]^d\to\R$ let
\[
  \norm{g}_{H^\zeta}
  =\max_{|a|\le\underline\zeta}\norm{\partial^ag}_\infty
  +\max_{|a|=\underline\zeta}\sup_{x\ne y}
  \frac{|\partial^ag(x)-\partial^ag(y)|}{\norm{x-y}^{\theta_\zeta}},
\]
so that for integer $\zeta$ the top derivative is Lipschitz, and let $H^\zeta$ be the set of $g$
with $\norm{g}_{H^\zeta}<\infty$.
Here $a\in\mathbb N_0^d$ is a multi-index, $|a|=\sum_j a_j$,
$a!=\prod_j a_j!$, and $\partial^a$ denotes the corresponding partial derivative.
The norm $\norm{\cdot}$ on $\R^d$ is Euclidean, and $\norm{\cdot}_\infty$ is the supremum norm.
Write $H^\zeta(L)=\{g:\norm{g}_{H^\zeta}\le L\}$ for the ball of radius $L$.

\paragraph{Model.}
We observe independent and identically distributed (i.i.d.) copies of $Z=(X,A,Y)$ with $X\in[0,1]^d$, $A\in\{0,1\}$, and
\[
  Y=\mu_0(X)+A\tau(X)+\varepsilon,\qquad \E(\varepsilon\mid X,A)=0,
  \qquad \pi(x)=\Prob(A=1\mid X=x).
\]
The covariate $X$ has density $f$. We assume $\pi\in H^\alpha$,
$\mu_0\in H^\beta$, and $\tau\in H^\gamma$, where $\alpha,\beta,\gamma>0$.
The bounds on the H\"older norms and the outcome moments are specified
for each model class below.
Define
\[
s=\frac{\alpha+\beta}{2},
\qquad
s_1=\frac{(\alpha\wedge 1)+(\beta\wedge 1)}{2},
\qquad
\bar\alpha=\max\{\alpha,\beta\},
\]
and
\[
\rho_{\mathrm{TS}}
=
\frac{2(s+s_1)}
{d+4s_1+\frac{2ds_1}{\gamma}},
\qquad
  \rho_{\mathrm S}=\frac{1}{1+d/(2\gamma)+d/(4s)}.
\]
The low-smoothness regime is $s<d\gamma/(4\gamma+2d)$.
The subscript $\mathrm{TS}$ denotes the two-scale construction used for rough design.

\subsection{Two design settings: rough and uniform}
\label{sec:designs}

\paragraph{The rough-design class.}
For rough design, fix density bounds
$0<\underline f\le1\le\overline f<\infty$ and assume
$\underline f\le f\le\overline f$ almost everywhere on $[0,1]^d$.
No smoothness is imposed on $f$.
Let $\cP_{\rm rough}(\alpha,\beta,\gamma)$ denote the laws in the model above satisfying
this density bound and
\[
  \begin{gathered}
  \epsilon\le\pi\le1-\epsilon,\qquad
  \pi\in H^\alpha(L_\pi),\quad \mu_0\in H^\beta(L_\mu),\quad
  \tau\in H^\gamma(L_\tau),\\
  \mathop{\rm ess\,sup}_{x,a}\E(Y^2\mid X=x,A=a)\le M_2.
  \end{gathered}
\]
Here $\epsilon\in(0,1/2)$ is the overlap constant, which bounds treatment probabilities
away from zero and one. The constants $\underline f,\overline f,L_\pi,L_\mu,L_\tau,M_2$
are positive and fixed as $n$ varies; their dependence is suppressed in the class notation.

\paragraph{The uniform-design setting.}
We also consider the uniform design, with $f\equiv1$ on $[0,1]^d$.
The estimator is constructed without using the uniformity assumption.
\section{Rough design}
\label{sec:rough}
We work over $\cP_{\rm rough}(\alpha,\beta,\gamma)$. We first construct an
estimator and state its risk bound, then give the corresponding minimax lower bound.
\subsection{Upper bound}
\label{sec:rough-upper}
\paragraph{Stencil construction.}
Fix an integer $q>\alpha+\beta+\gamma+1$ and set
$D_q=\binom{d+q}{q}$ and $m_{\rm st}=D_q+1$.
Choose $D_q$ points $z_1=0,z_3,\ldots,z_{m_{\rm st}}$ at which evaluation
uniquely determines every polynomial of total degree at most $q$.
Choose sufficiently small neighborhoods $U_i$ of $z_i$,
$3\le i\le m_{\rm st}$, so that this property holds for every
selection of one point from each box together with the origin.

Split the $n$ observations $Z_i=(X_i,A_i,Y_i)$ into
$m_{\rm st}$ disjoint groups $I_1,\ldots,I_{m_{\rm st}}$ of approximately equal size.
For a kernel $K$ of one observation from each group, define the product-sample average
\begin{align*}
    \mathbb{U}_{n,m_{\rm st}}K=\frac{1}{\prod_{j=1}^{m_{\rm st}}|I_j|}\sum_{i_j\in I_j,1\leq j\leq m_{\rm st}}K(Z_{i_1},...,Z_{i_{m_{\rm st}}})
\end{align*}
Write $B(x,t)=[0,1]^d\cap\{u:\norm{u-x}_\infty\le t\}$.
Within this subsection, $\widehat Q$, $\varphi_h$, and the sample subsets $I_j$
refer to the rough-design estimator.

\begin{samepage}
\noindent An ordered tuple $(x_1,\ldots,x_{m_{\rm st}})$ belongs to the stencil event
$\mathcal C_{h,\ell,r}(x_0)$ if
\begin{align*}
    x_1&\in B(x_0,h)\\
    x_2&\in B(x_1,\ell)\\
    x_i&\in x_1+rU_i \quad 3\leq i\leq m_{\rm st}
\end{align*}
\end{samepage}
For an accepted tuple, let $V(x_{1:m_{\rm st}})$ be the $D_q \times m_{\rm st}$ polynomial evaluation matrix. Its rows are indexed by monomials of total degree at most $q$, evaluated at $(x_j-x_1)/r$. The matrix has rank $D_q$. Define the polynomial map as
\[
\varphi_h(x)
=
\left(
\frac{1}{a!}
\left(
\frac{x-x_0}{h}
\right)^a
\right)_{a\in\mathcal A_q}
\in \mathbb{R}^{D_q}
\]
where $\mathcal A_q=\{a\in\mathbb N_0^d:|a|\le q\}$.
Order the coordinates with $a=0$ first, and let $e_0\in\R^{D_q}$ select this coordinate.
Let
\[
w=(w_1,\ldots,w_{m_{\rm st}})^\top
\]
be the unit null vector of $V(x_{1:m_{\rm st}})$, with the sign fixed by $w_1>0$.

\paragraph{Estimator.}
Define
\begin{align*}
    &C_A=\sum_{j=1}^{m_{\rm st}}w_jA_j \quad C_Y=\sum_{j=1}^{m_{\rm st}} w_jY_j\\
    &C_{A\varphi}=\sum_{j=1}^{m_{\rm st}} w_jA_j\varphi_h(X_j) \quad \overline{\varphi}_w=\sum_{j=1}^{m_{\rm st}} w_j^2 \varphi_h(X_j)
\end{align*}
Write $\mathcal C=\mathcal C_{h,\ell,r}(x_0)$. Define
\begin{align*}
    \hat{b}&=\mathbb{U}_{n,m_{\rm st}}[1_\mathcal{C}\overline{\varphi}_wC_AC_Y]\\
    \hat{Q}&= \mathbb{U}_{n,m_{\rm st}}[1_\mathcal{C}\overline{\varphi}_wC_AC_{A\varphi}^\top]
\end{align*}
Let $s_{\min}$ denote the smallest singular value.
Let $q_0=h^d\ell^dr^{d(m_{\rm st}-2)}$, and let $\kappa_Q>0$ be a fixed constant such that
$s_{\min}(\E_P\widehat Q)\ge\kappa_Qq_0$ uniformly over the model class for the bandwidths used below.
The final estimator is
\[
\widetilde{\theta}
=
\begin{cases}
\widehat Q^{-1}\widehat b,
& s_{\min}(\widehat Q)\ge \kappa_Q q_0/2,\\[4pt]
0,
& \text{otherwise},
\end{cases}
\qquad
\widehat{\tau}_{h,\ell,r}(x_0)
=
\operatorname{clip}_{[-L_\tau,L_\tau]}
\left(e_0^\top\widetilde{\theta}\right).
\]
Here $\operatorname{clip}_{[a,b]}(t)=\max\{a,\min\{t,b\}\}$, and
$\norm{\tau}_\infty\le L_\tau$ by the model-class definition.

\paragraph{Risk bound.}
The following theorem gives the risk bound for this estimator.

\begin{theorem}[Rough-design upper bound]
\label{thm:rough-upper}
For suitable bandwidths $h,\ell,r$ and all sufficiently large $n$, the estimator
defined above satisfies
\[
\sup_{P\in\cP_{\rm rough}(\alpha,\beta,\gamma)}
\E_P|\widehat\tau_{h,\ell,r}(x_0)-\tau_P(x_0)|
\lesssim\begin{cases}
n^{-\rho_{\rm TS}}, & s<\dfrac{d\gamma}{4\gamma+2d},\\[6pt]
n^{-\gamma/(2\gamma+d)}, & s\ge\dfrac{d\gamma}{4\gamma+2d}.
\end{cases}
\]
\end{theorem}

\subsection{Lower bound}
\label{sec:rough-lower}
We next give a lower bound on the worst-case risk of any estimator.

\begin{theorem}[Rough-design lower bound]
\label{thm:rough-lower}
For every $\omega>0$ and all sufficiently large $n$,
\[
\inf_{\widehat\tau}\sup_{P\in\cP_{\rm rough}(\alpha,\beta,\gamma)}
\E_P|\widehat\tau(x_0)-\tau_P(x_0)|\gtrsim
\begin{cases}
n^{-(\rho_{\rm TS}+\omega)}, & s<\dfrac{d\gamma}{4\gamma+2d},\\[6pt]
n^{-\gamma/(2\gamma+d)}, & s\ge\dfrac{d\gamma}{4\gamma+2d}.
\end{cases}
\]
The infimum is over all measurable estimators based on $n$ i.i.d. observations.
\end{theorem}

\section{Uniform design}
\label{sec:smooth}

We first derive a conditional upper bound for general design densities
and then verify its conditions under the uniform design.
The lower bound is established under the uniform design.

\subsection{Upper bound}
\label{sec:smooth-upper}

We describe a second-order local polynomial R-learner \citep{kennedy2024minimax} for $\tau(x_0)$ at a fixed
point $x_0\in[0,1]^d$. The design distribution and the nuisance functions are estimated on
training data $\mathcal T$ independent of the estimation sample. The results are stated for a
general design density; the uniform design is the special case $f\equiv1$.

\paragraph{Assumptions on the distribution.}
Assume $\epsilon\le\pi\le1-\epsilon$, $\pi\in H^\alpha$, $\mu_0\in H^\beta$, $\tau\in H^\gamma$
with $\norm{\pi}_{H^\alpha}+\norm{\mu_0}_{H^\beta}+\norm{\tau}_{H^\gamma}\le L_{\rm mod}$,
$\norm{\mu_0}_\infty+\norm{\tau}_\infty\le M$, and
\[
  \E\{\exp(t\varepsilon)\mid X,A\}\le\exp(\sigma^2t^2/2),\qquad t\in\R.
\]
The marginal law $F$ of $X$ has a density $f$ with $0<\underline f\le f\le\overline f<\infty$ on $[0,1]^d$.
Define $\cP_{\mathrm{unif}}(\alpha,\beta,\gamma;\epsilon,L,M)$ as the laws in the model
of \cref{sec:model} with $f\equiv1$, $\epsilon\le\pi\le1-\epsilon$,
$\norm{\pi}_{H^\alpha}+\norm{\mu_0}_{H^\beta}+\norm{\tau}_{H^\gamma}\le L$,
$\norm{\mu_0}_\infty+\norm{\tau}_\infty\le3M$, and $|\varepsilon|\le3M$ almost surely,
where $L,M>0$ are fixed. Every law in this class satisfies
these assumptions with $\underline f=\overline f=1$, $L_{\rm mod}=L$, and with $M$ and $\sigma$ both replaced
by $3M$.

\paragraph{Kernel and local polynomial basis.}
Let $K_h(x_0,u)=h^{-d}2^{-d}\one\{\norm{u-x_0}_\infty\le h\}$,
$r_\gamma=\lfloor\gamma\rfloor$, $\mathcal A_{r_\gamma}=\{a\in\mathbb N_0^d:|a|\le r_\gamma\}$,
$p=|\mathcal A_{r_\gamma}|$, and
\[
  \varphi_{h,x_0}(u)
  =\Big(\frac1{a!}\prod_{j=1}^d\Big(\frac{u_j-x_{0,j}}h\Big)^{a_j}\Big)_{a\in\mathcal A_{r_\gamma}}\in\R^p,
  \qquad
  \varphi(t)=(t^a/a!)_{a\in\mathcal A_{r_\gamma}},
\]
with the coordinate $a=0$ first and $e_0\in\R^p$ the corresponding unit vector. Let
$B_{x_0}(h)=[0,1]^d\cap\prod_j[x_{0,j}-h,x_{0,j}+h]$. For $\varsigma\in\{-1,1\}^d$ let
$I_\varsigma=\prod_j I_{\varsigma_j}$ with $I_{1}=[0,1/2]$ and $I_{-1}=[-1/2,0]$, and set
\[
  \lambda_*=\min_{\varsigma\in\{-1,1\}^d}\lambda_{\min}\Big(\int_{I_\varsigma}\varphi(t)\varphi(t)^\top\,dt\Big)>0,
  \qquad
  c_Q=2^{-d}\underline f\,\epsilon(1-\epsilon)\lambda_* .
\]

\paragraph{Local mesh and estimated projection kernel.}
Let $0<\ell\le h\le1$. Partition each coordinate interval of $B_{x_0}(h)$ into intervals of
length in $[\ell/2,\ell]$ and take the induced partition of $B_{x_0}(h)$ into cells. Fix
$q_{\rm poly}\ge\max\{\lfloor\alpha\rfloor,\lfloor\beta\rfloor\}$, let $W_{x_0,\ell}$ be the space
of functions supported on $B_{x_0}(h)$ that are polynomials of degree at most $q_{\rm poly}$ on
each cell, and let
\[
  V_{x_0,\ell}=W_{x_0,\ell}+\sum_{a\in\mathcal A_{r_\gamma}}\varphi_{h,x_0,a}W_{x_0,\ell}.
\]
Let $\widehat F$ be a finite nonnegative measure on $[0,1]^d$ built from $\mathcal T$, and
$\widehat\mu_{x_0}(du)=K_h(x_0,u)\,d\widehat F(u)$. Fix a spanning set $b$ of $V_{x_0,\ell}$,
remove redundant vectors with respect to $L_2(\widehat\mu_{x_0})$ by a fixed rule, and write $b$
for the kept vectors. Define
\[
  \widehat\Omega_{x_0}=\int b(u)b(u)^\top\,d\widehat\mu_{x_0}(u),
  \qquad
  \widehat P_{x_0}(u,v)=b(u)^\top\widehat\Omega_{x_0}^{-1}b(v).
\]
Only the restriction of $\widehat F$ to $B_{x_0}(h)$ enters this projection;
$\widehat F$ need not have total mass one.

\paragraph{Nuisance estimators and the estimator.}
Let $\widehat\pi$ and $\widehat\mu_0$ be estimators built from $\mathcal T$, and write
$\Delta_\pi=\pi-\widehat\pi$ and $\Delta_\mu=\mu_0-\widehat\mu_0$.
In this upper-bound construction, $n$ denotes the size of the estimation sample, excluding training data;
$I_1,I_2$ and $\widehat Q$ are defined for the estimator in this section.
Split the estimation sample into independent parts $I_1,I_2$ with
$\underline c_In\le|I_j|\le\overline c_In$. For $i\in I_1$, $j\in I_2$, write
$\varphi_i=\varphi_{h,x_0}(X_i)$, $K_i=K_h(x_0,X_i)$, and likewise for $j$. Let
\begin{align*}
  \widehat Q&=\frac{1}{|I_1|}\sum_{i\in I_1}K_i\varphi_i\varphi_i^\top A_i\{A_i-\widehat\pi(X_i)\}\\
  &\quad-\frac{1}{|I_1||I_2|}\sum_{i\in I_1}\sum_{j\in I_2}
  K_i\varphi_i\{A_i-\widehat\pi(X_i)\}K_j\widehat P_{x_0}(X_i,X_j)A_j\varphi_j^\top,\\
  \widehat R&=\frac{1}{|I_1|}\sum_{i\in I_1}K_i\varphi_i\{A_i-\widehat\pi(X_i)\}\{Y_i-\widehat\mu_0(X_i)\}\\
  &\quad-\frac{1}{|I_1||I_2|}\sum_{i\in I_1}\sum_{j\in I_2}
  K_i\varphi_i\{A_i-\widehat\pi(X_i)\}K_j\widehat P_{x_0}(X_i,X_j)\{Y_j-\widehat\mu_0(X_j)\}.
\end{align*}
Fix $M_{\rm clip}\ge\norm{\tau}_\infty$ uniformly over the model class.
We take $M_{\rm clip}=M$ under the general distribution assumptions and
$M_{\rm clip}=3M$ for $\cP_{\mathrm{unif}}(\alpha,\beta,\gamma;\epsilon,L,M)$.
Let $G=\{s_{\min}(\widehat Q)\ge c_Q/2\}$,
\[
  \widehat\theta=
  \begin{cases}
  \widehat Q^{-1}\widehat R,& G,\\
  0,& G^c,
  \end{cases}
  \qquad
  \widehat\tau_F(x_0)=e_0^\top\widehat\theta,
  \qquad
  \widetilde\tau_F(x_0)=\max\{-M_{\rm clip},\min\{\widehat\tau_F(x_0),M_{\rm clip}\}\}.
\]

\paragraph{Constants.}
Let
\[
  \mathfrak C=
  \{d,\alpha,\beta,\gamma,\epsilon,\underline f,\overline f,L_{\rm mod},L_{\rm nuis},M,M_{\rm clip},\sigma,
  \underline c_I,\overline c_I,q_{\rm poly}\},
\]
where $L_{\rm nuis}$ is the constant in condition (2) below. Constants implicit in $\lesssim$
depend only on $\mathfrak C$ and on $C_{\rm F2}$ in condition (F2) below, and not on
$n,h,\ell,x_0,P$, or the training sample. For $g$ on $[0,1]^d$ and $x\in[0,1]^d$ write
$\norm{g}_{x,h}^2=\int K_h(x,u)g(u)^2\,dF(u)$.

\begin{theorem}[Conditional risk bound]
\label{thm:smooth-cond}
Fix $x_0\in[0,1]^d$ and $0<\ell\le h\le1$. Suppose the following hold uniformly over the
distribution class.
\begin{enumerate}
\item[(1)] The assumptions on the distribution above hold with the stated constants.
\item[(2)] Conditionally on $\mathcal T$,
$\epsilon/2\le\widehat\pi\le1-\epsilon/2$ and $\norm{\Delta_\pi}_{H^\alpha}+\norm{\Delta_\mu}_{H^\beta}\le L_{\rm nuis}$.
\item[(3)] $\ell^\alpha$ is small enough that $C_{\rm stab}\ell^\alpha\le c_Q/8$, where
$C_{\rm stab}=C_{\rm stab}(\mathfrak C)$ is a fixed positive stability constant depending only on $\mathfrak C$.
\item[(4)] Conditionally on $\mathcal T$, $\widehat F\ll F$, and with $r_F=d\widehat F/dF$
and $\Delta_F=\norm{r_F-1}_\infty$,
\begin{enumerate}
\item[(F1)] $\Delta_F\le c_F$, where $c_F=c_F(\mathfrak C)\in(0,1/2]$ satisfies
$C_{\rm stab}c_F\le c_Q/8$;
\item[(F2)] with $a_\pi=\norm{\Delta_\pi}_{x_0,h}$ and
$a_\mu=\norm{\Delta_\mu}_{x_0,h}$,
\[
  \Delta_Fa_\pi(a_\mu+h^\gamma)\le C_{\rm F2}(h^\gamma+\ell^{2s}).
\]
\end{enumerate}
\end{enumerate}
Let $\E^{\rm fin}$ denote expectation over the estimation sample only, with
$(\widehat F,\widehat\pi,\widehat\mu_0)$ fixed. Then there is a constant $C$, depending only
on $\mathfrak C$ and $C_{\rm F2}$, such that
\[
  \E^{\rm fin}|\widehat\tau_F(x_0)-\tau(x_0)|^2
  \le C\Big(h^{2\gamma}+\ell^{4s}+\frac{1}{nh^d}+\frac{1}{n^2h^d\ell^d}\Big).
\]
Consequently, with $h^\gamma=n^{-\rho_{\rm S}}$ and $\ell=h^{\gamma/(2s)}$, in the low-smoothness
regime,
\[
  \E^{\rm fin}|\widehat\tau_F(x_0)-\tau(x_0)|^2\lesssim n^{-2\rho_{\rm S}}.
\]
Both bounds are uniform in $x_0\in[0,1]^d$, in $P$, and in the training realization.
\end{theorem}

The proof uses conditions (2) and (4) only on $B_{x_0}(h)$. So \cref{thm:smooth-cond} remains
true if the H\"older norms in (2) and the supremum in $\Delta_F$ are taken over any set
containing $B_{x_0}(h)$ instead of $[0,1]^d$.

\begin{remark}[Known design distribution]
If $F$ is known, take $\widehat F=F$. Then $\Delta_F=0$ and (F1), (F2) hold trivially.
\end{remark}

\paragraph{Unconditional bound.}
Define the event
\[
  \mathcal E_{\rm tr}
  =\left\{
  \widehat F\ll F,\ \text{(F1) and (F2) hold},\
  \epsilon/2\le\widehat\pi\le1-\epsilon/2,\
  \norm{\Delta_\pi}_{H^\alpha}+\norm{\Delta_\mu}_{H^\beta}\le L_{\rm nuis}
  \right\},
\]
which is measurable with respect to $\mathcal T$, and assume that uniformly over $P$,
\begin{equation}
  \Prob_P(\mathcal E_{\rm tr}^c)\le\delta_{\rm tr}
  \label{eq:smooth-U1}
\end{equation}
for a deterministic $\delta_{\rm tr}$, which may depend on the sizes of the training samples.

\begin{theorem}[Unconditional risk bound]
\label{thm:smooth-uncond}
Suppose conditions (1) and (3) of \cref{thm:smooth-cond} hold uniformly over the model class
and \eqref{eq:smooth-U1} holds. Let
$R_{n,h,\ell}=h^{2\gamma}+\ell^{4s}+(nh^d)^{-1}+(n^2h^d\ell^d)^{-1}$. Then, with the constant
$C$ of \cref{thm:smooth-cond},
\[
  \Prob_P\left[\E_P\{|\widetilde\tau_F(x_0)-\tau(x_0)|^2\mid\mathcal T\}\le CR_{n,h,\ell}\right]
  \ge1-\delta_{\rm tr},
\]
and
\[
  \E_P|\widetilde\tau_F(x_0)-\tau(x_0)|^2\le CR_{n,h,\ell}+4M_{\rm clip}^2\delta_{\rm tr}.
\]
Both bounds are uniform in $x_0\in[0,1]^d$ and in $P$.
\end{theorem}

\paragraph{Uniform but unknown design.}
We now show that conditions (2) and (4) of \cref{thm:smooth-cond} hold with high probability
when the design is uniform, for estimators that do not use this fact. Fix $x_0\in(0,1)^d$ and
$r_0>0$ with $\operatorname{dist}_\infty(x_0,\partial[0,1]^d)\ge2r_0$, and let
$B_0=x_0+[-r_0,r_0]^d$. For $g$ on $B_0$, $\norm{g}_{H^\zeta(B_0)}$ is the H\"older norm
of Section~\ref{sec:model} with the suprema taken over $B_0$. The training data consist of
two independent samples, $\mathcal T=(\mathcal T_N,\mathcal D_f)$, where
$\mathcal T_N=\{(X_k^{\rm tr},A_k^{\rm tr},Y_k^{\rm tr}):k\le N\}$ is used for $\widehat\pi$ and
$\widehat\mu_0$ and $\mathcal D_f=\{X_k^f:k\le m_f\}$ is a covariate sample used for
$\widehat F$. Fix a kernel $K_{\rm tr}\in C^{J+1}(\R^d)$,
$J=\lceil\bar\alpha\rceil$, with $K_{\rm tr}\ge0$, $\int K_{\rm tr}=1$, and support in
$[-1,1]^d$, and write $K_\lambda(t)=\lambda^{-d}K_{\rm tr}(t/\lambda)$. Define
\[
  \widehat f(u)=\frac1{m_f}\sum_{k=1}^{m_f}K_{\lambda_f}(u-X_k^f),
  \qquad
  S_\lambda[\xi](x)=\frac1N\sum_{k=1}^NK_\lambda(x-X_k^{\rm tr})\,\xi(Z_k^{\rm tr})
\]
for bounded functions $\xi$ of $Z=(X,A,Y)$, let
$d\widehat F(u)=\one\{u\in[0,1]^d\}\widehat f(u)\,du$ without renormalization, and
let
\[
  \widehat\pi(x)=\max\Big\{\frac\epsilon2,\min\Big\{\frac{S_\lambda[A](x)}{S_\lambda[1](x)},1-\frac\epsilon2\Big\}\Big\},
  \qquad
  \widehat\mu_0(x)=\frac{S_{\lambda_\mu}[(1-A)Y](x)}{S_{\lambda_\mu}[1-A](x)},
\]
with $\widehat\pi=\epsilon/2$ where $S_\lambda[1]=0$, $\widehat\mu_0=0$ where its denominator
is zero, and $\lambda_\mu=\lambda$ if $\alpha\ge\beta$ and $\lambda_\mu=r_0$ if
$\alpha<\beta$.

\begin{theorem}[Verification under uniform design]
\label{thm:smooth-unif}
Let $P\in\cP_{\mathrm{unif}}(\alpha,\beta,\gamma;\epsilon,L,M)$, let $0<\ell\le h\le r_0$,
and let $0<\delta\le1/2$. Take
\[
  \lambda=N^{-\eta}\ \text{with }0<\eta<\frac1{d+2\bar\alpha},
  \qquad
  \lambda_f\in[1/m_f,r_0].
\]
There are constants $C_1,C_2,L_{\rm nuis}^*$, depending only
on $(d,\alpha,\beta,\gamma,\epsilon,L,M,q_{\rm poly},K_{\rm tr},r_0)$, such that if
\[
  N^{1-\eta(d+2\bar\alpha)}\ge C_1\log(N/\delta),
  \qquad
  m_f\lambda_f^d\ge C_2\,h^{-2\gamma}\log(m_f/\delta),
\]
then with probability at least $1-\delta$ over $\mathcal T$, the following
hold with all norms and suprema taken over $B_0$:
\begin{enumerate}
\item[(a)] $\epsilon/2\le\widehat\pi\le1-\epsilon/2$ and
$\norm{\Delta_\pi}_{H^\alpha(B_0)}+\norm{\Delta_\mu}_{H^\beta(B_0)}\le L_{\rm nuis}^*$;
\item[(b)] $\widehat F\ll F$, $\Delta_F=\sup_{B_0}|\widehat f-1|\le c_F$, and
$\Delta_F\,a_\pi(a_\mu+h^\gamma)\le h^\gamma$.
\end{enumerate}
That is, conditions (2) and (4) of \cref{thm:smooth-cond} hold on $B_0$ with
$L_{\rm nuis}=L_{\rm nuis}^*$ and $C_{\rm F2}=1$. The probability bound is uniform over
$P\in\cP_{\mathrm{unif}}(\alpha,\beta,\gamma;\epsilon,L,M)$.
\end{theorem}

Since the design density is one, on $B_0$ the mean of $S_\lambda[\xi]$ is the convolution
$K_\lambda*g$ with $g(u)=\E(\xi\mid X=u)$. Convolution does not increase H\"older norms, so the
H\"older norms of the ratios defining $\widehat\pi$ and $\widehat\mu_0$ are controlled by those
of $\pi$ and $\mu_0$, and the stochastic parts have H\"older norm of order
$\lambda^{-\bar\alpha}\sqrt{\log(N/\delta)/(N\lambda^d)}$, which the condition on $\eta$ makes
small. Likewise $\widehat f$ is unbiased on $B_0$, so $\Delta_F$ is of order
$\sqrt{\log(m_f/\delta)/(m_f\lambda_f^d)}$. This is where uniformity of the design is used, and
it is why $x_0$ is taken in the interior.

\paragraph{Risk bound.}

Take $h^\gamma=n^{-\rho_{\rm S}}$, $\ell=h^{\gamma/(2s)}$, and
$\delta=n^{-2\rho_{\rm S}}$ in \cref{thm:smooth-unif}, with $N=m_f=n$ and $\lambda_f=r_0$. In the
low-smoothness regime $\rho_{\rm S}<1/2$, so both sample size conditions of \cref{thm:smooth-unif}
hold for large $n$. Then \cref{thm:smooth-uncond}, with $\delta_{\rm tr}=n^{-2\rho_{\rm S}}$ and with the
conditions localized to $B_0$ as noted after \cref{thm:smooth-cond}, gives the following result.

\begin{corollary}[Uniform-design upper bound]
\label[corollary]{cor:smooth-upper}
Suppose $s<d\gamma/(4\gamma+2d)$ and $x_0\in(0,1)^d$.
With the bandwidths and training sample sizes chosen above, for all sufficiently large $n$,
\[
  \sup_{P\in\cP_{\mathrm{unif}}(\alpha,\beta,\gamma;\epsilon,L,M)}
  \E_P|\widetilde\tau_F(x_0)-\tau(x_0)|^2\le Cn^{-2\rho_{\rm S}}.
\]
The estimator does not use the design density, and $M_{\rm clip}=3M$.
\end{corollary}

The total sample size is $n_{\rm tot}=n+N+m_f=3n$;
replacing $n$ by $n_{\rm tot}$ changes only the constant in the rate.

\subsection{Lower bound}
\label{sec:smooth-lower}

We next establish a lower bound under uniform design, using a binary-outcome subclass.
Define this subclass by
\[
  \cP^{\mathrm{bin}}_{\mathrm{unif}}(\alpha,\beta,\gamma;\epsilon)
  =\{P\in\cP_{\mathrm{unif}}(\alpha,\beta,\gamma;\epsilon,1,1):
    P(Y\in\{0,1\})=1\}.
\]

\begin{theorem}[Uniform-design lower bound]
\label{thm:smooth-lower}
Let $d\ge1$, let $\alpha,\beta,\gamma>0$ satisfy $s<d\gamma/(4\gamma+2d)$, let
$\epsilon\in(0,1/2)$, and let $x_0\in(0,1)^d$. There is a constant $c>0$, depending only on
$(d,\alpha,\beta,\gamma,\epsilon,x_0)$, such that for all $n$,
\[
  \inf_{\widehat\tau}\ \sup_{P\in\cP^{\mathrm{bin}}_{\mathrm{unif}}(\alpha,\beta,\gamma;\epsilon)}
  \E_P\big|\widehat\tau(x_0)-\tau_P(x_0)\big|\ \ge\ c\,n^{-\rho_{\rm S}},
\]
where the infimum is over all measurable functions of $n$ i.i.d. observations.
\end{theorem}

Since
\[
  \cP^{\mathrm{bin}}_{\mathrm{unif}}(\alpha,\beta,\gamma;\epsilon)
  \subseteq\cP_{\mathrm{unif}}(\alpha,\beta,\gamma;\epsilon,L,M)
  \qquad\text{for all }L\ge1,\ M\ge1,
\]
the same lower bound holds over $\cP_{\mathrm{unif}}(\alpha,\beta,\gamma;\epsilon,L,M)$.

\cref{thm:smooth-lower} shows that no estimator can do better than $n^{-\rho_{\rm S}}$ in mean
absolute error. Together with \cref{cor:smooth-upper}, this shows that
in the low-smoothness regime $s<d\gamma/(4\gamma+2d)$ the minimax rate for
$\tau(x_0)$ at an interior point under uniform design is $n^{-\rho_{\rm S}}$ with
$\rho_{\rm S}=1/\{1+d/(2\gamma)+d/(4s)\}$.

\section{Discussion}
Our results show that the minimax rate for estimating the CATE at a point depends on the covariate design.

Future work could characterize the minimax rate when $f\in H^\nu(L_f)$ for
fixed $\nu,L_f>0$, where $\nu$ measures design-density smoothness. This would
describe the intermediate setting between rough and smooth designs. Another
direction is to extend the analysis to average treatment effect (ATE) estimation
and other causal estimands.
\label{sec:discussion}

\section*{Proofs}
Lean 4 code for the proofs is available on GitHub at
\url{https://github.com/houyx22/causal_lowerbound}.
Parts of the proofs were developed with AI assistance.

\nocite{*}

\bibliographystyle{plainnat}
\bibliography{refs}

@article{rosenbaum1983central,
  author  = {Rosenbaum, Paul R. and Rubin, Donald B.},
  title   = {The central role of the propensity score in observational studies for causal effects},
  journal = {Biometrika},
  year    = {1983},
  volume  = {70},
  number  = {1},
  pages   = {41--55},
  doi     = {10.1093/biomet/70.1.41}
}

@article{heckman1998matching,
  author  = {Heckman, James J. and Ichimura, Hidehiko and Todd, Petra},
  title   = {Matching as an econometric evaluation estimator},
  journal = {The Review of Economic Studies},
  year    = {1998},
  volume  = {65},
  number  = {2},
  pages   = {261--294},
  doi     = {10.1111/1467-937X.00044}
}

@article{abadie2006large,
  author  = {Abadie, Alberto and Imbens, Guido W.},
  title   = {Large sample properties of matching estimators for average treatment effects},
  journal = {Econometrica},
  year    = {2006},
  volume  = {74},
  number  = {1},
  pages   = {235--267},
  doi     = {10.1111/j.1468-0262.2006.00655.x}
}

@article{abadie2011bias,
  author  = {Abadie, Alberto and Imbens, Guido W.},
  title   = {Bias-corrected matching estimators for average treatment effects},
  journal = {Journal of Business \& Economic Statistics},
  year    = {2011},
  volume  = {29},
  number  = {1},
  pages   = {1--11},
  doi     = {10.1198/jbes.2009.07333}
}

@inproceedings{gao2020minimax,
  author    = {Gao, Zijun and Han, Yanjun},
  title     = {Minimax optimal nonparametric estimation of heterogeneous treatment effects},
  booktitle = {Advances in Neural Information Processing Systems},
  year      = {2020},
  volume    = {33},
  pages     = {21751--21762},
  publisher = {Curran Associates, Inc.},
  url       = {https://proceedings.neurips.cc/paper/2020/hash/f75b757d3459c3e93e98ddab7b903938-Abstract.html}
}

@article{lin2023estimation,
  author  = {Lin, Zhexiao and Ding, Peng and Han, Fang},
  title   = {Estimation based on nearest neighbor matching: From density ratio to average treatment effect},
  journal = {Econometrica},
  year    = {2023},
  volume  = {91},
  number  = {6},
  pages   = {2187--2217},
  doi     = {10.3982/ECTA20598}
}

@article{robins1994estimation,
  author  = {Robins, James M. and Rotnitzky, Andrea and Zhao, Lue Ping},
  title   = {Estimation of regression coefficients when some regressors are not always observed},
  journal = {Journal of the American Statistical Association},
  year    = {1994},
  volume  = {89},
  number  = {427},
  pages   = {846--866},
  doi     = {10.1080/01621459.1994.10476818}
}

@article{bang2005doubly,
  author  = {Bang, Heejung and Robins, James M.},
  title   = {Doubly robust estimation in missing data and causal inference models},
  journal = {Biometrics},
  year    = {2005},
  volume  = {61},
  number  = {4},
  pages   = {962--973},
  doi     = {10.1111/j.1541-0420.2005.00377.x}
}

@article{lee2017doubly,
  author  = {Lee, Sokbae and Okui, Ryo and Whang, Yoon-Jae},
  title   = {Doubly robust uniform confidence band for the conditional average treatment effect function},
  journal = {Journal of Applied Econometrics},
  year    = {2017},
  volume  = {32},
  number  = {7},
  pages   = {1207--1225},
  doi     = {10.1002/jae.2574}
}

@article{nie2021quasi,
  author  = {Nie, Xinkun and Wager, Stefan},
  title   = {Quasi-oracle estimation of heterogeneous treatment effects},
  journal = {Biometrika},
  year    = {2021},
  volume  = {108},
  number  = {2},
  pages   = {299--319},
  doi     = {10.1093/biomet/asaa076}
}

@article{kennedy2023towards,
  author  = {Kennedy, Edward H.},
  title   = {Towards optimal doubly robust estimation of heterogeneous causal effects},
  journal = {Electronic Journal of Statistics},
  year    = {2023},
  volume  = {17},
  number  = {2},
  pages   = {3008--3049},
  doi     = {10.1214/23-EJS2157}
}

@article{kennedy2024minimax,
  author  = {Kennedy, Edward H. and Balakrishnan, Sivaraman and Robins, James M. and Wasserman, Larry},
  title   = {Minimax rates for heterogeneous causal effect estimation},
  journal = {The Annals of Statistics},
  year    = {2024},
  volume  = {52},
  number  = {2},
  pages   = {793--816},
  doi     = {10.1214/24-AOS2369}
}


\end{document}